\documentstyle[12pt,amsfonts]{article}
\newtheorem{theorem}{Theorem}[section]
\newtheorem{corollary}[theorem]{Corollary}
\newtheorem{definition}{Definition}
\newtheorem{example}[theorem]{Example}
\newtheorem{remark}[theorem]{Remark}

\newtheorem{proposition}[theorem]{Proposition}
\title{G.$\Lambda_s$-sets and G.$V_s$-sets\thanks{1991 Math.\
Subject Classification ---Primary: 54D30, 54A05; Secondary:
54H05, 54G99. \protect\newline Keywords and phrases --- generalized
closed sets, semi-open, semi-T$_{\frac{1}{2}}$space,
semi-$T_1$-space, semi-$R_0$-space. \protect\newline Research supported
partially by the Ella and Georg Ehrnrooth Foundation at Merita
Bank, Finland.}}
\author{Miguel Caldas Cueva\\Departamento de Matematica
Aplicada\\Universidade Federal Fluminense\\IMUFF-Rua Mario Santos
Braga\\s/n, CEP: 24020-140, Niteroi\\RJ-Brasil \and Julian
Dontchev\\Department of Mathematics\\University of Helsinki\\PL
4, Yliopistonkatu 15\\00014 Helsinki 10\\Finland}
\date{}
\begin{document}
\baselineskip=20pt plus 1pt minus 1pt
\maketitle
\begin{abstract}
In this paper we define the concepts of $g.\Lambda_s$-sets and
$g.V_s$-sets and we use them in order to obtain new
characterizations of semi-$T_1$-, semi-$R_0$- and
semi-T$_{\frac{1}{2}}$-spaces.
\end{abstract}

\section{Introduction}\label{s0}

Separation axioms stand among the most common and to a certain
extent the most important and interesting concepts in Topology.
One of the most well-known low separation axiom is the one which
requires that singletons are closed, i.e.\ $T_1$. In most
studies, spaces under consideration are `by default' $T_1$.

In Digital Topology \cite{KKM1} several spaces that fail to be
$T_1$ are important in the study of the geometric and topological
properties of digital images \cite{KR1,KK1,K1}. Such is the case
with the major building block of the digital n-space -- the {\em
digital line} or the so called {\em Khalimsky line}. This is the
set of the integers, $\mathbb Z$, equipped with the
topology $\cal K$, generated by ${\cal G}_{\cal K} = \{ \{ 2n-1,
2n, 2n+1 \} \colon n \in {\mathbb Z} \}$.

Although the digital line is neither a $T_1$-space nor an
$R_0$-space, it satisfies a couple of separation axioms which are
a bit weaker than $T_1$ and $R_0$, that is, the digital line is
both a semi-$T_1$-space and a semi-$R_0$-space. This inclines to
indicate that further knowledge of the behavior of topological
spaces satisfying these two weak separation axioms (and some
related ones) is required. This is indeed the intention of the
present paper.

\section{Preliminaries}\label{s1}

The concept of a semi-open set in a topological space was
introduced by N.~Levine in 1963 \cite{L1}. If $(X,\tau)$ is a
topological space and $A \subseteq X$, then $A$ is {\em
semi-open} \cite{L1} if there exists $O \in \tau$ such that $O
\subseteq A \subseteq {\rm Cl}(O)$, where ${\rm Cl}(O)$ denotes
closure of $O$ in $(X,\tau)$. The complement $A^c$ of a semi-open
set $A$, is called {\em semi-closed} and the {\em semi-closure}
of a set $A$ denoted by ${\rm sCl}(A)$, is the intersection of
all semi-closed sets containing $A$. 

The separation axioms $R_0$ was introduced by Davis in \cite{D2}.
It requires that every open set contains the closures of its
points. In 1975, Maheshwari and Prasad \cite{MP2} introduced the
class of semi-$R_0$-spaces studied later by Di Maio \cite{D1} and
by Jankovi\'{c} and Reilly \cite{JR1}. A topological space
$(X,\tau)$ is called a {\em semi-$R_0$-space} if every semi-open
set contains the semi-closure of each of its singletons. One can
easily observe that a space $(X,\tau)$ is a semi-$R_0$-space if
and only if every semi-open set is union of semi-closed sets. In
this paper we focus our attention precisely on the sets which are
union of semi-closed sets and study their basic properties.
Although that the first impression might be that the separation
axiom semi-$R_0$ is rather weak, one needs to consider the fact
that $T_4$-spaces, even ultraconnected spaces, need not be
semi-$R_0$ (the easiest example is probably a Sierpinski space).

A generalized class of closed sets was considered by Maki in 1986
\cite{M1}. He investigated the sets that can be represented as
union of closed sets and called them {\em $V$-sets}. Complements
of $V$-sets, i.e., sets that are intersection of open sets are
called {\em $\Lambda$-sets} \cite{M1}. In connection to
semi-$R_0$-spaces, observe that $R_0$-spaces are precisely the
spaces where the closed sets form a network for the topology,
i.e., the spaces where every open set is a $V$-set. Every
$R_0$-space is a semi-$R_0$-space \cite{JR1} but not vice versa.

The family of all semi-open (resp.\ semi-closed) sets in
$(X,\tau)$ will be denoted by $SO(X,\tau)$ (resp.\ $SC(X,\tau)$).
In this note, we introduce and characterize the concepts of
$\Lambda_s$-set, $V_s$-set, $g.\Lambda_s$-set and $g.V_s$-set in
a topological space $(X,\tau)$. In this article we give new
characterizations of semi-$T_{\frac{1}{2}}$-, semi-$R_0$- and
semi-$T_1$-spaces in terms $V_s$- and $g.V_s$-sets.

\section{$\Lambda_s$-sets and $V_s$-sets}\label{s2}

\begin{definition}\label{d1}
{\em Let $B$ be a subset of a topological space $(X,\tau)$. We
define the subsets $B^{\Lambda_s}$ and $B^{V_s}$ as follows:

$B^{\Lambda_s} = \bigcap \{O/O\supseteq B$, $O\in SO(X,\tau )\}$
and $B^{V_s} = \bigcup \{F/F\subseteq B$, $F^c \in SO(X,\tau)\}$.

In \cite{DM1,MP2}, $B^{\Lambda_s}$ is called the {\em
semi-kernel} of $B$.}
\end{definition}

\begin{proposition}\label{p1}
Let $A,B$ and $\{ B_\lambda \colon \lambda \in \Omega\}$ be
subsets of a topological space $(X,\tau)$. Then the following
properties are valid:

(a) $B \subseteq B^{\Lambda_s}$;

(b) If $A \subseteq B$, then $A^{\Lambda_s} \subseteq B^{\Lambda
_s}$;

(c) $B_{}^{\Lambda_s\Lambda_s} = B^{\Lambda_s}$;

(d) $[\bigcup \limits_{\lambda \in
\Omega}^{}B_\lambda]^{\Lambda_s} =
\bigcup \limits_{\lambda \in \Omega} B_\lambda^{\Lambda_s}$;

(e) If $A\in SO(X,\tau)$, then $A = A^{\Lambda_s}$;

(f) $(B^c)^{\Lambda_s}=(B^{V_s})^c$;

(g) $B^{V_s} \subseteq B$;

(h) If $B \in SC(X,\tau)$, then $B=B^{V_s}$;

(i) $[\bigcap\limits_{\lambda \in \Omega }B_\lambda ]^{\Lambda
_s} \subseteq \bigcap\limits_{\lambda \in \Omega }B_\lambda
^{\Lambda _s}$;

(j) $[\bigcup\limits_{\lambda \in \Omega }B_\lambda ]^{V_s}
\supseteq \bigcup\limits_{\lambda \in \Omega }B_\lambda ^{V_s}$
\end{proposition}

{\em Proof.} (a) Clear by Definition~\ref{d1}.

(b) Suppose that $x\notin B^{\Lambda _s}$. Then there exists a
subset $O \in SO(X,\tau)$ such that $O \supseteq B$ with $x
\not\in O$. Since $B \supseteq A$, then $x \not\in A^{\Lambda
_s}$ and thus $A^{\Lambda _s} \subseteq B^{\Lambda_s}$.

(c) Follows from (a) and Definition~\ref{d1}.

(d) Suppose that there exists a point $x$ such that $x \not\in
[\bigcup\limits_{\lambda \in \Omega }B_\lambda ]^{\Lambda _s}.$
Then, there exists a subset $O\in SO(X,\tau )$ such that
$\bigcup\limits_{\lambda \in \Omega }B_\lambda \subseteq O$ and
$x\notin O.$ Thus, for each $\lambda \in \Omega $ we have $x
\not\in B_\lambda^{\Lambda _s}$. This implies that $x \not\in
\bigcup\limits_{\lambda \in \Omega }B_\lambda ^{\Lambda _s}$.
Conversely, suppose that there exists a point $x \in X$ such that
$x \not\in \bigcup\limits_{\lambda \in \Omega }B_\lambda^{\Lambda
_s}.$ Then by Definition~\ref{d1}, there exist subsets $O_\lambda
\in SO(X,\tau)$ (for all $\lambda \in \Omega$) such that $x
\not\in O_\lambda$, $B_\lambda \subseteq O_\lambda $. Let $O =
\bigcup\limits_{\lambda \in \Omega }O_\lambda $. Then we have
that $x \not\in \bigcup\limits_{\lambda \in \Omega }O_\lambda$,
$\bigcup\limits_{\lambda \in \Omega }B_\lambda \subseteq O$ and
$O\in SO(X,\tau)$. This implies that $x \not\in
[\bigcup\limits_{\lambda \in \Omega }B_\lambda ]^{\Lambda _s}$.
Thus, the proof of (d) is completed.

(e) By Definition~\ref{d1} and since $A \in SO(X,\tau)$, we have
$A^{\Lambda_s} \subseteq A$. By (a) we have that $A^{\Lambda _s}
= A$. 

(f) $(B^{V_s})^c=\bigcap \{F^c/F^c\supseteq B^c, F^c\in
SO(X,\tau)\} = (B^c)^{\Lambda _s}$.

(g) Clear by Definition~\ref{d1}.

(h) If $B \in SC(X,\tau)$, then $B^c \in SO(X,\tau )$. By (e) and
(f): $B^c=(B^c)^{\Lambda _s}=(B^{V_s})^c$. Hence $B=B^{V_s}$.

(i) Suppose that there exists a point $x$ such that $x \not\in
\bigcap\limits_{\lambda \in \Omega }B_\lambda^{\Lambda _s}$.
Then, there exists $\lambda \in \Omega $ such that $x \not\in
B_\lambda^{\Lambda_s}$. Hence there exists $\lambda \in \Omega$
and $O \in SO(X,\tau)$ such that $O \supseteq B_\lambda $ and $x
\not\in O$. Thus $x \not\in [\bigcap\limits_{\lambda \in \Omega
}B_\lambda ]^{\Lambda_s}$.

(j) $[\bigcup\limits_{\lambda \in \Omega }B_\lambda
]^{V_s}=[((\bigcup\limits_{\lambda \in \Omega }B_\lambda
)^c)^{\Lambda _s}]^c=[(\bigcap\limits_{\lambda \in \Omega
}B_\lambda ^c)^{\Lambda _s}]^c\supseteq [\bigcap\limits_{\lambda
\in \Omega }(B_\lambda ^c)^{\Lambda _s}]^c = [\bigcap
\limits_{\lambda \in \Omega }(B_\lambda ^{V_s})^c]^c$
$=\bigcup\limits_{\lambda \in \Omega }B_\lambda ^{V_s}$ (by (f)
and (h)). $\Box$

\begin{remark}\label{r1}
{\em In general $(B_1\bigcap B_2)^{\Lambda _s}\neq B_1^{\Lambda
_s} \bigcap B_2^{\Lambda _s}$, as the following example shows.}
\end{remark}

\begin{example}\label{e1}
{\em Let $(X,\tau)$ be as in (\cite[Example 2.9]{M1}) i.e., let
$X=\{a,b,c\}$ and $\tau = \{\emptyset, \{a\}, \protect\newline \{b,c\},
X\}$. Let $B_1=\{b\}$ and $B_2=\{c\}.$ Then, $(B_1 \bigcap
B_2)^{\Lambda _s}=\emptyset$ but $B_1^{\Lambda_s}\bigcap
B_2^{\Lambda
_s} = \{b,c\}.$}
\end{example}

\begin{definition}\label{d2}
{\em In a topological space $(X,\tau)$, a subset $B$ is a
{\em $\Lambda_s$-set} (resp.\ {\em $V_s$-set}) of $(X,\tau)$ if
$B=B^{\Lambda_s}$ (resp.\ $B=B^{V_s})$.}
\end{definition}

\begin{remark}\label{r2}
{\em By Proposition~\ref{p1} (e) and (h) we have that:

(a) If $B$ is a $\Lambda$-set or if $B \in SO(X,\tau)$, then $B$
is a $\Lambda_s$-set.

(b) If $B$ is a $V$-set or if $B\in SC(X,\tau )$, then $B$ is a
$V_s$-set.}
\end{remark}

\begin{proposition}\label{p2}
(a) The subsets $\emptyset$ and $X$ are $\Lambda_s$-sets and
$V_s$-sets.

(b) Every union of $\Lambda_s$-sets (resp.\ $V_s$-sets) is a
$\Lambda_s$-set (resp.\ $V_s$-set).

(c) Every intersection of $\Lambda_s$-sets (resp.\ $V_s$-sets)
is a $\Lambda_s$-set (resp.\ $V_s$-set).

(d) A subset $B$ is a $\Lambda_s$-set if and only if $B^c$ is a
$V_s$-set.
\end{proposition}

{\em Proof.} (a) and (d) are obvious.

(b) Let $\{B_\lambda \colon \lambda \in \Omega \}$ be a family
of $\Lambda_s$-set in a topological space $(X,\tau)$. Then by
Definition~\ref{d2} and Proposition~\ref{p1} (d), $\bigcup
\limits_{\lambda \in \Omega} B_\lambda = \bigcup \limits_{\lambda
\in \Omega} B_\lambda^{\Lambda_s} = [\bigcup\limits_{\lambda \in
\Omega }B_\lambda ]^{\Lambda _s}$.

(c) Let $\{B_\lambda :\lambda \in \Omega \}$ be a family of
$\Lambda_s$-set in $(X,\tau)$. Then by Proposition~\ref{p1} (h)
and Definition~\ref{d2} $[\bigcap\limits_{\lambda \in \Omega
}B_\lambda ]^{\Lambda _s}\subseteq \bigcap\limits_{\lambda \in
\Omega }B_\lambda ^{\Lambda_s}=\bigcap\limits_{\lambda \in \Omega
}B_\lambda$. Hence by Proposition~\ref{p1} (a) $\bigcap
\limits_{\lambda \in \Omega }B_\lambda = [\bigcap\limits_{\lambda
\in \Omega }B_\lambda ]^{\Lambda_s}$. $\Box$

Recall that a space topological $(X,\tau)$ is called a {\em
semi-$T_1$-space} \cite{MP1} if to each pair of distinct points
$x,y$ of $(X,\tau)$ there corresponds a semi-open set $A$
containing $x$ but not $y$ and a semi-open set $B$ containing $y$
but not $x$, or equivalently, $(X,\tau)$ is a semi-$T_1$-space
if and only if every singleton is semi-closed.

\begin{example}
{\em The digital line $({\mathbb Z},{\cal K})$ is an example of
a semi-$T_1$ space and a semi-$R_0$-space which is neither $T_1$
nor $R_0$. Since all even singletons are closed, they are
trivially semi-closed. On the other hand the odd integers are
regular open (but not closed) and hence semi-closed too. Thus,
$({\mathbb Z},{\cal K})$ is semi-$T_1$ and every semi-open sets
is the union of all of its semi-closed singletons. So, the
digital line is a semi-$R_0$-space. On the other hand the
isolated points in the digital line (i.e., the odd integers) can
not be expressed as union of closed sets, which implies that the
digital line is not an $R_0$-space.}
\end{example}

\begin{proposition}\label{p3}
A topological space $(X,\tau)$ is a semi-$T_1$-space if and only
if every subset is a $\Lambda_s$-set (or equivalently a
$V_s$-set).
\end{proposition}

{\em Proof.} Let $B$ be a subset of a semi-$T_1$-space
$(X,\tau)$. Suppose that there exists a point $x\in X$ such that
$x \not\in B$. Then, $\{x\}^c$ is a semi-open set containing $B$.
Then, by Definition~\ref{d1} $B^{\Lambda_s}\subseteq \{x\}^c$.
This implies $x \not\in B^{\Lambda_s}$. Hence we have $B^{\Lambda
_s}\subseteq B$ and $B^{\Lambda_s} = B$ (Proposition~\ref{p1}
(a)). 

For the converse, if $x \in X$, then $X \setminus \{x\}^c$ is due
to assumption a $\Lambda_s$-set. Hence, its complement $\{ x \}$
is union of semi-closed sets and thus semi-closed. This shows
that $X$ is a semi-$T_1$-space. $\Box$

\begin{corollary}
Every semi-$T_1$-space is a semi-$R_0$-space.
\end{corollary}

{\em Proof.} The definition of semi-$R_0$-spaces requires that
every semi-open set is a $V_s$-set. $\Box$

\begin{example}\label{e3a}
{\em Since indiscrete spaces (with at least two points) are
semi-$R_0$, then the separation axiom semi-$R_0$ is strictly
below semi-$T_1$. However, it is interesting to mention that even
closed subspaces of semi-$T_1$-spaces need not be semi-$R_0$ and
that semi-$T_1$-spaces need not be $R_0$. Let $X = \{ a,b,c \}$
and $\tau = \{ \emptyset, \{ a \}, \{ b \}, \{ a,b \}, X \}$.
Observe that $X$ is semi-$T_1$ and that the closed subspace $A
= \{ a,c \}$ is not a semi-$R_0$-space. Moreover, $X$ is not
$R_0$.}
\end{example}

Recall that a subset $A$ is called {\em simply-open} if $A$ is
union of an open and a nowhere dense set. If $A \subseteq {\rm
Int} ({\rm Cl} (A))$, then $A$ is called {\em locally dense} or
{\em preopen}. Sets which are dense in some regular closed
subspace are called {\em semi-preopen} or {\em $\beta$-sets}.

\begin{theorem}
For a topological space $(X,\tau)$ the following conditions are 
equivalent:

(1) $X$ is a semi-$T_1$-space;

(2) Every locally dense (= preopen) subspace is a $V_s$-set;

(3) Every $\beta$-open (= semi-preopen) subspace is a
$V_s$-set.
\end{theorem}

{\em Proof.} (1) $\Rightarrow$ (3) Follows from
Proposition~\ref{p3}. 

(3) $\Rightarrow$ (2) Obvious, since every locally dense set is
$\beta$-open.

(2) $\Rightarrow$ (1) Let $x \in X$. It is well-known that every
singleton is either locally dense or nowhere dense \cite{JR1}.
If $\{ x \}$ is locally dense, then by (2), $\{ x \}$ is union
of
semi-closed sets and hence semi-closed. If $\{ x \}$ is nowhere
dense, then it is clearly semi-closed. Thus every singleton of
$X$ is semi-closed and consequently $X$ is a semi-$T_1$-space.
$\Box$

\begin{theorem}
For a topological space $(X,\tau)$ the following conditions are
equivalent:

(1) $X$ is a semi-$R_0$-space;

(2) Every simply-open (= locally semi-closed) subspace is a 
$V_s$-set;

(3) Every open subspace is a $V_s$-set.
\end{theorem}

{\em Proof.} Since every open set is semi-open and since every
semi-open set is simply-open, (1) $\Rightarrow$ (3) and (2)
$\Rightarrow$ (1) are obvious.

(3) $\Rightarrow$ (2) If $A \subseteq X$ is simply-open, then $A
= U \bigcup N$, where $U \in \tau$ and $N$ is nowhere dense. By
(3), $U$ is a $V_s$-set. Since every nowhere dense set is
semi-closed, then by Proposition~\ref{p2} $A$ is a $V_s$-set.
$\Box$

\begin{example}\label{e3}
{\em In the Euclidean plane $({\mathbb R}^{2},\tau)$, every
singleton $\{x\}$ is a $\Lambda_s$-set by Proposition~\ref{p3}.
However, $\{x\}$ is not semi-open in $({\mathbb R}^{2},\tau)$.
Thus the converse of Proposition~\ref{p1} (e) is not true in
general.}
\end{example}

\section{G.$\Lambda_S$-sets and g.$V_S$-sets}\label{s3}

In this section, by using the $\Lambda_s$-operator and
$V_s$-operator, we introduce the classes of generalized
$\Lambda_s$-sets (= $g.\Lambda_s$-sets) and generalized
$V_s$-sets (= $g.V_s$-sets) as an analogy of the sets 
introduced by H.~Maki \cite{M1}.

\begin{definition}\label{d3}
{\em In a topological space $(X,\tau)$, a subset $B$ is called
a {\em $g.\Lambda_s$-set} of $(X,\tau)$ if $B^{\Lambda_s}
\subseteq F$ whenever $B \subseteq F$ and $F$ is semi-closed.}
\end{definition}

\begin{definition}\label{d4}
{\em In a topological space $(X,\tau)$, a subset $B$ is called
a {\em $g.V_s$-set} of $(X,\tau)$ if $B^c$\ is a
$g.\Lambda_s$-set of $(X,\tau)$.}
\end{definition}

\begin{remark}
{\em By $D^{\Lambda_s}$ (resp.\ $D^{V_s})$ we will denote the
family of all $g.\Lambda_s$-sets (resp.\ $g.V_s$-sets) of
$(X,\tau)$.}
\end{remark}

\begin{proposition}\label{p32}
Let $(X,\tau)$ be a topological space. Then: 

(a) Every $\Lambda_s$-set is a $g.\Lambda_s$-set. 

(b) Every $V_s$-set is a $g.V_s$-set. 

(c) If $B_\lambda \in D^{\Lambda _s}$ for all $\lambda \in
\Omega$, then $\bigcup \limits_{\lambda \in \Omega }B_\lambda \in
D^{\Lambda_s}$.

(d) If $B_\lambda \in D^{V_s}$ for all $\lambda \in \Omega$, then
$\bigcap\limits_{\lambda \in \Omega }B_\lambda \in D^{V_s}$. 
\end{proposition}

{\em Proof.} (a) Follows from Definition~\ref{d2} and
Definition~\ref{d3}.

(b) Let $B$ be a $V_s$-set subset of $X$. Then, $B = B^{V_s}$.
By Proposition~\ref{p1} (f) $(B^c)^{\Lambda_s} = (B^{V_s})^c =
B^c.$ Therefore, by (a) and Definition~\ref{d4}, $B$ is a
$g.V_s$-set.

(c) Let $B_\lambda \in D^{\Lambda _s}$ for all $\lambda \in
\Omega$. Then, by Proposition~\ref{p1} (d) $[\bigcup
\limits_{\lambda \in \Omega} B_\lambda ]^{\Lambda_s} = \bigcup
\limits_{\lambda \in \Omega }B_\lambda ^{\Lambda _s}$. Hence, by
hypothesis and Definition~\ref{d3}, $\bigcup \limits_{\lambda \in
\Omega }B_\lambda \in D^{\Lambda_s}$.

(d) Follows from (c) and Definition~\ref{d4}. $\Box$ 

In general the intersection of two $g.\Lambda_s$-sets is not a
$g.\Lambda_s$-sets as shown by the following example.

\begin{example}\label{e33}
{\em Let $X=\{a,b,c\}$ and $\tau = \{ \emptyset, \{a,b\}, X \}$.
If $A = \{a,c\}$ and $B = \{b,c\}$ (as in \cite[Example
3.3]{M1}). Then $A$ and $B$ are $g.\Lambda_s$-sets, but $A
\bigcap B = \{c\}$ is not a $g.\Lambda_s$-set. We have:
$D^{\Lambda_s} = \{\emptyset, \{a\}, \{b\}, \{a,b\}, \{a,c\},
\{b,c\}, X \}$ and $D^{V_s} = \{ \emptyset, \{a\}, \{b\}, \{c\},
\{a,c\}, \{b,c\}, X \}$.}
\end{example}

The following example shows that the converse of
Proposition~\ref{p32} (a) (resp.\ (b)) is not true in general.

\begin{example}\label{e34}
{\em Let $(X,\tau )$ be the space in Example~\ref{e33}. The
subset $A = \{a,c\}$ is a $g.\Lambda_s$-set but it is not a
$\Lambda_s$-set.}
\end{example}

\begin{remark}
{\em By Remark~\ref{r2} and Proposition~\ref{p32} we have that:

(i) If $A \in SO(X,\tau)$, then $A$ is a $g.\Lambda_s$-set;

(ii) If $A \in SC(X,\tau)$, then $A$ is a $g.V_s$-set.}
\end{remark}

\begin{proposition}\label{p36}
Let $(X,\tau )$ be a topological space.

(a) For each $x \in X$, $\{x\}$ is a semi-open set or $\{x\}^c$
is a $g.\Lambda_s$-set of $(X,\tau)$.

(b) For each $x \in X$, $\{x\}$ is a semi-open set or $\{x\}$ is
a $g.V_s$-set of $(X,\tau)$.
\end{proposition}

{\em Proof.} Suppose that $\{x\}$ is not semi-open. Then only
semi-closed set $F$ containing $\{x\}^c$ is $X$. Thus
$(\{x\}^c)^{\Lambda_s} \subseteq F = X$ and $\{x\}^c$ is a
$g.\Lambda_s$-set of $(X,\tau)$.

(b) Follows from (a) and Definition~\ref{d3}. $\Box$

\begin{corollary}
For a topological space $(X,\tau)$, the Cantor-Bendixson
derivative $D(X)$ is the set of all $g.V_s$-singletons of
$(X,\tau)$.
\end{corollary}

\begin{proposition}\label{p37}
If $A$ is a $g.\Lambda_s$-set of a topological space $(X,\tau)$
and $A \subseteq B \subseteq A^{\Lambda_s}$, then $B$ is a
$g.\Lambda_s$-set of $(X,\tau)$.
\end{proposition}

{\em Proof.} Since $A \subseteq B \subseteq A^{\Lambda _s}$, we
have $A^{\Lambda_s} = B^{\Lambda_s}$ by Proposition~\ref{p1} (b).
Let $F$ be any semi-closed subset of $(X,\tau)$ such that $B
\subseteq F$. Then, we have $B^{\Lambda_s} = A^{\Lambda_s}
\subseteq F$, since $A \subseteq B$ and $A$ is $g.\Lambda_s$-set.
$\Box$

In the following propositions we give a characterization of
$g.V_s$-sets (Definition~\ref{d4}) by using $V_s$-operations and
we obtain results concerning such subsets.

\begin{proposition}\label{p38}
A subset $B$ of a topological space $(X,\tau )$ is a $g.V_s$-set
if and only if $U \subseteq B^{V_s}$ whenever $U \subseteq B$ and
$U\in SO(X,\tau)$.
\end{proposition}

{\em Proof.} {\em Necessity.} Let $U$ be a semi-open subset of
$(X,\tau)$ such that $U \subseteq B$. Then since $U^c$ is
semi-closed and $U^c \supseteq B^c$, we have $U^c \supseteq
(B^c)^{\Lambda_s}$ by Definition~\ref{d3} and
Definition~\ref{d4}. Hence by Proposition~\ref{p1} (f)
$U^c\supseteq (B^{V_s})^c$. Thus, $U \subseteq B^{V_s}$.

{\em Sufficiency.} Let $F$ be a semi-closed subset of $(X,\tau)$
such that $B^c \subseteq F$. Since $F^c$ is semi-open and $F^c
\subseteq B$, by assumption we have $F^c \subseteq B^{V_s}$.
Then, $F \supseteq (B^{V_s})^c = (B^c)^{\Lambda_s}$ by
Proposition~\ref{p1} (f), and $B^c$ is a $g.\Lambda_s$-set, i.e.,
$B$ is a $g.V_s$-set. $\Box$

As consequence of Proposition~\ref{p38}, we have:

\begin{corollary}\label{p39}
Let $B$ be a $g.V_s$-set in a topological space $(X,\tau)$. Then,
for every semi-closed set $F$ such that $B^{V_s} \bigcup B^c
\subseteq F$, $F=X$ holds.
\end{corollary}

{\em Proof.} The assumption $B^{V_s} \bigcup B_{}^c \subseteq F$
implies $(B^{V_s})^c \bigcap B \supseteq F^c$. Since $B$ is a
$g.V_s$-set, then by Proposition~\ref{p38}, we have $B^{V_s}
\supseteq F^c$ and hence $(B^{V_s})^c\subseteq F$
and $\emptyset = (B^{V_s})^c \bigcap B^{V_s} \supseteq F^c$.
Therefore, we have $X=F$. $\Box$

\begin{corollary}\label{310}
Let $B$ be a $g.V_s$-set of $(X,\tau)$. Then $B^{V_s} \bigcup
B^c$ is a semi-closed set if and only if $B$ is a $V_s$-set.
\end{corollary}

{\em Proof.} {\em Necessity.} By Proposition~\ref{p39}, $B^{V_s}
\bigcup B^c=X$. Thus $(B^{V_s})^c \bigcap B = \emptyset$. Hence,
by Proposition~\ref{p1} (g) $B=B^{V_s}$. {\em Sufficiency} is
obvious. $\Box$

\begin{proposition}\label{p311}
Let $B$ be a subset of topological space $(X,\tau)$ such that
$B^{V_s}$ is semi-closed. If $X=F$ holds for every semi-closed
subset $F$ such that $F \supseteq B^{V_s} \bigcup B^c$, then $B$
is a $g.V_s$-set.
\end{proposition}

{\em Proof.} Let $U$ be a semi-open subset contained in $B$.
According to assumption, $B^{V_s} \bigcup U^c$ is semi-closed
such that $B^{V_s} \bigcup B^c \subseteq B^{V_s} \bigcup U^c$.
It follows that $B^{V_s} \bigcup U^c = X$ and hence $U
\subseteq B^{V_s}$. By Proposition~\ref{p38}, $B$ is a
$g.V_s$-set. $\Box$

\section{Characterization of semi-$T_{\frac{1}{2}}$
spaces}\label{s4}

After the work of N.~Levine \cite{L1} on semi-open sets, various
mathematicians turned their attention to the generalizations of
various concepts in topology by considering semi-open sets
instead of open sets. In this direction, P.~Bhattacharyya and
B.~K.~Lahiri \cite{BL1} defined the concept of semi-generalized
closed (= sg-closed) sets of a topological space in terms of
semi-open sets. In recent years, the class of
semi-T$_{\frac{1}{2}}$-spaces has been of some interest, (i.e.,
the spaces where the classes of semi-closed sets and the
sg-closed sets coincide). In this section we give a new
characterization of semi-T$_{\frac{1}{2}}$-spaces by using
$g.V_s$-sets. In order to achieve our purpose, we recall the
following definitions (see also \cite{BL1,C1,M1}).

\begin{definition}\label{d5}
{\em A subset $B$ of a topological space $(X,\tau)$ is said to
be {\em semi-generalized closed set} (written briefly as
$sg$-closed) \cite{BL1} if ${\rm sCl}(B) \subseteq O$ holds
whenever $B \subseteq O$ and $O\in SO(X,\tau)$. Every semi-closed
sets is sg-closed but the converse is not always true
\cite{BL1}.}
\end{definition}

\begin{definition}\label{d6}
{\em A topological space $(X,\tau)$ is said to be a {\em
semi-T$_{\frac{1}{2}}$-space} \cite{BL1} if every $sg$-closed set
in $(X,\tau)$ is semi-closed in $(X,\tau)$.}
\end{definition}

In the following theorem we give a characterization of the class
of semi-$T_{\frac{1}{2}}$ by using $g.V_s$-sets.

\begin{theorem}\label{t41}
Let $(X,\tau)$ be a topological space. Then the following
statements are equivalent:

(a) $(X,\tau)$ is a semi-T$_{\frac{1}{2}}$-space.

(b) Every $g.V_s$-set is a $V_s$-set.
\end{theorem}

{\em Proof.} $(a) \Rightarrow (b)$: Suppose that there exists a
$g.V_s$-set $B$ which is not a $V_s$-set. Since $B^{V_s}
\subseteq B$ $(B^{V_s}\neq B)$, then there exists a point $x \in
B$ such that $x \not\in B^{V_s}$. Then the singleton $\{x\}$ is
not semi-closed. According to Proposition 4.3 of \cite{SMB1},
$\{x\}^c$ is a $sg$-closed set. On the other hand, we have that
$\{x\}$ is not semi-open (since $B$ is a $g.V_s$-set, $x \not\in
B^{V_s}$ and Proposition~\ref{p38}). Therefore, we have that
$\{x\}^c$ is not semi-closed but it is a sg-closed set. This
contradicts to the assumption that $(X,\tau)$ is a
semi-T$_{\frac{1}{2}}$-space.

$(b) \Rightarrow (a)$: Suppose that $(X,\tau)$ is not a
semi-T$_{\frac{1}{2}}$-space. Then, there exists a $sg$-closed
set $B$ which is not semi-closed. Since $B$ is not semi-closed,
there exists a point $x$ such that $x \not\in B$ and $x \in {\rm
sCl}(B)$. By Proposition~\ref{p36}, we have the singleton $\{x\}$
is a semi-open set or it is a $g.V_s$-set. When $\{x\}$ is
semi-open, we have $\{x\} \bigcap B \neq \emptyset$ because $x
\in {\rm sCl}(B)$. This is a contradiction. Let us consider the
case: $\{x\}$ is a $g.V_s$-set. If $\{x\}$ is not semi-closed,
we have $\{x\}^{V_s} = \emptyset$ and hence $\{x\}$ is not a
$V_s$-set. This contradicts to $(b)$. Next, if $\{x\}$ is
semi-closed, we have $\{x\}^c \supseteq {\rm sCl}(B)$ (i.e., $x
\not\in {\rm sCl}(B)).$ In fact, the semi-open set $\{x\}^c$
contains the set $B$ which is a $sg$-closed set. Then, this also
contradicts to the fact that $x\in {\rm sCl}(B)$. Therefore
$(X,\tau)$ is a semi-T$_{\frac{1}{2}}$-space. $\Box$

\
E-mail: {\tt gmamccs@vm.uff.br}, {\tt dontchev@cc.helsinki.fi}
\
\end{document}